\documentclass[epsfig,latexsym,amsfonts,twoside]{article}
\usepackage{amssymb,amsmath,amscd,epsf,epsfig}

\def\part#1{\frac{\partial\phantom{#1}}{\partial#1}}
\newtheorem{thm}{Theorem}

\newtheorem{lem}[thm]{Lemma}

\newenvironment{prf}{\begin{trivlist}\item[]{\bf Proof} }%
{\hfill $\Box$ \end{trivlist}}
\newenvironment{dfn}{\begin{trivlist}\item[]{\bf Definition}\em }%
{\end{trivlist}}
\newenvironment{rmk}{\begin{trivlist}\item[]{\bf Remark} }%
{\end{trivlist}}
\newenvironment{exm}{\begin{trivlist}\item[]{\bf Example} }%
{\end{trivlist}}

\def\Z{\ifmmode{{\mathbb Z}}\else{${\mathbb Z}$}\fi}
\def\Q{\ifmmode{{\mathbb Q}}\else{${\mathbb Q}$}\fi}
\def\R{\ifmmode{{\mathbb R}}\else{${\mathbb R}$}\fi} 
\def\C{\ifmmode{{\mathbb C}}\else{${\mathbb C}$}\fi} 
\def\P{\ifmmode{{\mathbb P}}\else{${\mathbb P}$}\fi} 

\def\H{\ifmmode{{\mathrm H}}\else{${\mathrm H}$}\fi} 

\def\B{\ifmmode{{\cal B}}\else{${\cal B}$}\fi} 
\def\E{\ifmmode{{\cal E}}\else{${\cal E}$}\fi} 
\def\F{\ifmmode{{\cal F}}\else{${\cal F}$}\fi} 
\def\K{\ifmmode{{\cal K}}\else{${\cal K}$}\fi} 
\def\L{\ifmmode{{\cal L}}\else{${\cal L}$}\fi} 
\def\M{\ifmmode{{\cal M}}\else{${\cal M}$}\fi} 
\def\N{\ifmmode{{\cal N}}\else{${\cal N}$}\fi} 
\def\O{\ifmmode{{\cal O}}\else{${\cal O}$}\fi} 
\def\U{\ifmmode{{\cal U}}\else{${\cal U}$}\fi}
\def\X{\ifmmode{{\cal X}}\else{${\cal X}$}\fi} 

\def\Br{\ifmmode{{\mathrm{Br}}}\else{${\mathrm{Br}}$}\fi} 
\def\OG{\ifmmode{\widetilde{\cal M}_4}\else{$\widetilde{\cal M}_4$}\fi} 
\def\D{\ifmmode{{\cal{D}}_{\mathrm{coh}}^b}\else{${{\cal{D}}_{\mathrm{coh}}^b}$}\fi}
\def\Shah{\ifmmode{\amalg\hspace*{-3.5pt}\amalg}\else{$\amalg\hspace*{-3.5pt}\amalg$}\fi}

\newcommand{\normalization}{\raisebox{-20pt}{
                 \begin{picture}(60,60)(-30,-30)
                 \put(-30,30){\line(1,0){60}} 
                 \put(-30,30){\line(0,-1){60}} 
                 \put(30,30){\line(0,-1){60}} 
                 \put(-30,-30){\line(1,0){60}} 
                 \put(-20,27){\vector(1,-3){18}}
                 \put(0,27){\vector(1,-3){18}}
                 \end{picture}}}
 
\newcommand{\base}{\raisebox{0pt}{
                 \begin{picture}(60,10)(-30,-5)
                 \put(-30,0){\line(1,0){60}}
                 \end{picture}}}

\newcommand{\singularI}{\raisebox{-20pt}{
                 \begin{picture}(60,60)(-30,-20)
                 \put(-30,40){\line(1,0){60}} 
                 \put(-30,40){\line(0,-1){60}} 
                 \put(30,40){\line(0,-1){60}} 
                 \put(-30,-20){\line(1,0){60}} 
                 \put(-30,-20){\line(0,-1){20}}
                 \put(30,-20){\line(0,-1){20}}
                 \end{picture}}}

\newcommand{\singularII}{\raisebox{-20pt}{
                 \begin{picture}(60,60)(-30,-40)
                 \put(-30,20){\line(1,0){60}} 
                 \put(-30,20){\line(0,-1){60}} 
                 \put(30,20){\line(0,-1){60}} 
                 \put(-30,-40){\line(1,0){60}} 
                 \put(-30,20){\line(0,1){20}}
                 \put(30,20){\line(0,1){20}}
                 \end{picture}}}

\begin{document}

\title{On the discriminant locus of a Lagrangian
  fibration} 
\author{Justin Sawon}
\date{July, 2006}
\maketitle

\begin{abstract}
Let $X\rightarrow\P^n$ be an irreducible holomorphic symplectic
manifold of dimension $2n$ fibred over $\P^n$. Matsushita proved that
the generic fibre is a holomorphic Lagrangian abelian variety. In this
article we study the discriminant locus $\Delta\subset\P^n$
parametrizing singular fibres. Our main result is a formula for the
degree of $\Delta$, leading to bounds on the degree when $X$ is a
four-fold.
\end{abstract}

\section{Introduction}

Due primarily to the work of Matsushita~\cite{matsushita99,
  matsushita00i}, much is now known about the structure of fibrations
on irreducible holomorphic symplectic manifolds. In particular, the
generic fibre must be a holomorphic Lagrangian complex torus and it is
expected that the base must be projective space. In fact, ten years
earlier Mukai~\cite{mukai88} already posed the question: when is a
fibration $X\rightarrow\P^n$ by $n$-dimensional complex tori a
holomorphic symplectic manifold? Our goal in this article is to find
restrictions on the degree of the discriminant locus
$\Delta\subset\P^n$ in the case that $X$ is holomorphic symplectic.

To begin with, we assume the fibration is the relative Jacobian of a
family of curves, where the curves degenerate in a controlled manner
over a generic point of $\Delta$ (they acquire a single node). We
prove that the degree of $\Delta$ is given by
$$\mathrm{deg}\Delta=24\left(n!\sqrt{\hat{A}}[X]\right)^{\frac{1}{n}}$$
where $\sqrt{\hat{A}}[X]$ the characteristic number of $X$ coming from
the square root of the $\hat{A}$-polynomial. This same characteristic
number arose in earlier work of Hitchin and the author~\cite{hs01},
and it appears to play a fundamental r{\^o}le in holomorphic
symplectic geometry. The important point in this case is that we have
a model for the singular fibre $X_t$ above a generic point
$t\in\Delta$. Indeed the above formula for $\mathrm{deg}\Delta$
readily generalizes to
projective $X$ fibred by principally polarized abelian varieties,
provided the singular fibre $X_t$ for generic $t\in\Delta$ conforms to
this same model (Theorem~\ref{principal}). We then generalize our
formula to fibrations by non-principally polarized abelian varieties,
whose singular fibres conform to slightly different models
(Theorem~\ref{nonprincipal}). All of these models come from toroidal
compactifications of moduli spaces of abelian varieties, due to
Igusa~\cite{igusa67} and Mumford~\cite{mumford75} (see
also Ash et al.~\cite{amrt75}). In addition to constructing these
compactifications, Mumford~\cite{mumford72} described degenerations of
abelian varieties which sit above the boundary. Our assumption is that
for generic $t\in\Delta$, the singular fibre $X_t$ is a (semi-stable)
degeneration of an abelian variety which occurs over a generic point
in a codimension one component of the boundary. We will of course give
an explicit description of how these degenerate abelian varieties
look.

It is worth noting that in four dimensions
Matsushita~\cite{matsushita01} has classified all possible singular
fibres that can occur over a generic point of $\Delta$, and given
local models. Some of these look like products of smooth and singular
elliptic curves, up to {\'e}tale cover, and occur in examples like the
Hilbert scheme $\mathrm{Hilb}^2S$ of two points on an elliptic K3
surface $S\rightarrow\P^1$ (Example 3.5 in~\cite{sawon03}). Excluding
such Lagrangian fibrations, where the generic fibre is a product of
elliptic curves, the singular fibres considered in this article are
the only ones known to occur in global examples. It would be good to
extend our results to allow any of the singular fibres on Matsushita's
list, though the existence of non-reduced components creates some
difficulties.

In~\cite{guan01} Guan proved that the characteristic numbers of a
holomorphic symplectic four-fold are bounded; so when $X$ is a
four-fold our formulae give bounds on the degree of $\Delta$. We briefly
indicate why such bounds might be useful. Suppose $X\rightarrow\P^2$
is a fibration by abelian surfaces with polarization of type
$(1,d)$. This leads to a morphism
$$\phi:\P^2\backslash\Delta\longrightarrow\mathcal{A}^{\circ}(1,d)$$
where $\mathcal{A}^{\circ}(1,d)$ is the moduli space of abelian
surfaces with this polarization. If the singular fibres $X_t$ for
generic $t\in\Delta$ are well-behaved, then this map can be extended
to a morphism between (partial) compactifications
$$\phi^*:\P^2\backslash\Delta_{\mathrm{sing}}\longrightarrow\mathcal{A}^{\circ
  *}(1,d).$$
The construction and properties of $\mathcal{A}^{\circ *}(1,d)$ when
$d$ is prime are well-described in the book by Hulek, Kahn, and
Weintraub~\cite{hkw93}. The hope then is that the degree of $\Delta$
can be used to control the degree of the morphism $\phi^*$, implying
finiteness of the number of deformation classes of holomorphic
symplectic four-folds which admit Lagrangian fibrations (c.f.\ the
comments at the end of the introduction in Todorov~\cite{todorov03}).

We do not pursue this direction in this paper. Instead we use Guan's
bounds to show that both $d$ and the degree of $\Delta$ can take a
limited number of values (Theorem~\ref{values}).

The author would like to thank Manfred Lehn, Ivan Smith, and Richard
Thomas for useful conversations, and the Max-Planck-Institut f{\"u}r
Mathematik (Bonn) for hospitality. This work was also supported by NSF
grant number 0305865.

\section{Good singular fibres}

In this article a {\em Lagrangian fibration} shall mean an irreducible
holomorphic symplectic manifold $X$ of dimension $2n$ which is fibred
over projective space
$$f:X\rightarrow\P^n.$$
Let $\Delta\subset\P^n$ be the discriminant locus over which the
Jacobian of $f$ drops rank; it is a divisor parametrizing singular
fibres of $f$. The singular locus $\Delta_{\mathrm{sing}}$ of $\Delta$
will be codimension at least two in $\P^n$, which will mean that it
can effectively be ignored in most of our calculations. We write 
$\Delta_{\mathrm{sm}}:=\Delta\backslash\Delta_{\mathrm{sing}}$ for the
smooth locus of $\Delta$.

Matsushita~\cite{matsushita99, matsushita00i} proved that the generic
fibre of $f$ must be a (holomorphic Lagrangian) complex torus. We
begin by describing an example where the fibres are Jacobians of genus
$n$ curves.

\begin{exm}(The Beauville-Mukai integrable system~\cite{beauville99})
Let $S$ be a K3 surface which contains a smooth genus $n$ curve $C$,
and assume for simplicity that the Picard group of $S$ is generated
(over $\mathbb{Z}$) by this curve. Then $C$ moves in an
$n$-dimensional linear system $|C|\cong\P^n$ and every curve in this
family
$$\mathcal{C}\rightarrow\P^n$$
is integral (reduced and irreducible). The relative compactified
Jacobian $X=\bar{J}_0(\mathcal{C}/\P^n)$ is then a (smooth) Lagrangian
fibration over $\P^n$. Here the compactified Jacobian $\bar{J}_0C_t$
of an integral curve $C_t$ is defined to be the moduli space of
rank-one torsion-free sheaves of Euler characteristic zero, i.e.\
degree $n-1$ (see D'Souza~\cite{dsouza79}).
\end{exm}

There are two features of this fibration to which we wish to draw
attention. Firstly, each smooth fibre contains a canonical theta
divisor $\Theta$ (the image of
$$\mathrm{Sym}^{n-1}C_t\rightarrow\mathrm{Pic}^{n-1}C_t=J_0C_t,$$
which can be defined without reference to a basepoint). When
$t\in\Delta_{\mathrm{sm}}$, the curve $C_t$ acquires a single node. In this
case too there is a (generalized) theta divisor $\Theta$ on
$\bar{J}_0C_t$ (for example, see Esteves~\cite{esteves97}). So we have
a relative theta divisor over $\P^n\backslash\Delta_{\mathrm{sing}}$,
whose closure gives a divisor $Y$ in $X$.
 
Secondly, consider the structure of a singular fibre $\bar{J}_0C_t$
for $t\in\Delta_{\mathrm{sm}}$. The following description of the
compactified Jacobian of a curve $C_t$ with a single node is well
known (see Igusa~\cite{igusa56};
Example (1) on page 83 of Oda and Seshadri~\cite{os79} describes the
genus two case, which can easily be generalized). Let $\tilde{C_t}$ be
the normalization of $C_t$. The normalization of $\bar{J}_0C_t$ is
then a certain $\P^1$-bundle over $J_0\tilde{C_t}$. The zero and
infinity sections $s_0$ and $s_{\infty}$ of the $\P^1-$bundle are
canonically isomorphic to $J_0\tilde{C_t}$, but we instead identify
them using a certain translation in $J_0\tilde{C_t}$. Then
$\bar{J}_0C_t$ is given by taking the $\P^1$-bundle and gluing $s_0$ and
$s_{\infty}$ using the above identification. 
$$\begin{array}{cccc}
 & s_{\infty} & & \\
\mathbb{P}^1 & \hookrightarrow & \normalization & \bar{J}_0C_t
\\
 & s_0 & \downarrow & \\
 & & \base & J_0\tilde{C_t} \\
\end{array}$$

\begin{dfn}
Let $X\rightarrow\P^n$ be a Lagrangian fibration by principally
polarized abelian varieties such that the generic singular fibre
$X_t$ for $t\in\Delta_{\mathrm{sm}}$ is obtained by gluing together
the zero and infinity sections of a $\P^1$-bundle over a principally
polarized abelian variety of dimension $n-1$, just as in the example
above. Then we say $X\rightarrow\P^n$ has {\em good} singular fibres.
\end{dfn}

\begin{rmk}
As mentioned in the introduction, Igusa~\cite{igusa67} and
Mumford~\cite{mumford75} constructed compactifications of the moduli
space of abelian varieties. Although this involves some choices, in
the principally polarized case there is just one boundary component of
codimension one. Moreover, Mumford~\cite{mumford72} also gave a
construction of degenerating abelian varieties; a generic point of the
boundary then corresponds to a degenerate abelian variety as described
above, i.e.\ a good singular fibre, which can therefore be regarded as
the generic semi-stable degeneration of a principally polarized
abelian variety.
\end{rmk}

For a Lagrangian fibration with good singular fibres we arrive at the
following picture of the local structure of the fibration
$f:X\rightarrow\P^n$ over $\Delta_{\mathrm{sm}}$. In a neighbourhood
of the singular locus of a fibre over $\Delta_{\mathrm{sm}}$ there
exist local coordinates $(z_1,\ldots,z_n,w_1,\ldots,w_n)$ on $X$ such
that $f$ is given by
$$f:(z_1,\ldots,z_n,w_1,\ldots,w_n)\mapsto (z_1w_1,z_2,\ldots,z_n).$$
Here $\Delta_{\mathrm{sm}}$ is given by the vanishing of the first
component, locally on $\P^n$.



\section{The Beauville-Bogomolov quadratic form}

Let $X$ be an irreducible holomorphic symplectic manifold of dimension
$2n$. There is a quadratic form $q_X$ on $\H^2(X,\Z)$ known as the
Beauville-Bogomolov quadratic form (see~\cite{beauville83}). This form
generalizes the intersection pairing on a K3 surface. We begin with
some formulae involving $q_X$, which may be found in Huybrechts' notes
in~\cite{ghj02}, for instance.

The Fujiki formula states that
\begin{eqnarray}
\label{fujiki1}
q_X(\alpha)^n & = & \mbox{const.}\int_X\alpha^{2n}
\end{eqnarray}
for all $\alpha\in\H^2(X,\Z)$, where the constant depends only on
$X$. Fujiki also proved that if $\eta\in\H^{4j}(X,\R)$ is
of pure Hodge type $(2j,2j)$ on $X$ and on all small deformations of
$X$ then 
$$q_X(\alpha)^{n-j}=\mbox{const.}\int_X\eta\alpha^{2(n-j)}$$
for all $\alpha\in\H^2(X,\Z)$, where the constant depends only on
$\eta$. In particular, the second Chern class $c_2(T_X)$ satisfies the
hypothesis and thus
\begin{eqnarray}
\label{fujiki2}
q_X(\alpha)^{n-1} & = & \mbox{const.}\int_Xc_2\alpha^{2n-2}.
\end{eqnarray}
Writing out Equations~(\ref{fujiki1}) and~(\ref{fujiki2}) for $\alpha$
and $\beta\in\H^2(X,\Z)$, we can eliminate $q_X(\alpha)$,
$q_X(\beta)$, and both constants to obtain
\begin{eqnarray}
\label{alpha_beta}
\left(\int_X\alpha^{2n}\right)^{n-1}\left(\int_Xc_2\beta^{2n-2}\right)^n
 & = &
 \left(\int_X\beta^{2n}\right)^{n-1}\left(\int_Xc_2\alpha^{2n-2}\right)^n.
\end{eqnarray}
This equation will eventually yield a formula for the degree of the
discriminant locus.

We return to the situation of the previous section. Thus we
have a Lagrangian fibration $f:X\rightarrow\P^n$ with a divisor $Y$
which restricts to the theta divisor on each smooth fibre and to the
generalized theta divisor on a generic singular fibre (over
$\Delta_{\mathrm{sm}}$). There is also a divisor $L$ given by pulling
back a hyperplane from $\P^n$. We denote the holomorphic symplectic
form by $\sigma$. Substituting $\alpha=\sigma+t_1\bar{\sigma}$ and
$\beta=Y+t_2L$ into Equation~(\ref{alpha_beta}), and then
comparing coefficients of $(t_1t_2)^{n(n-1)}$ gives
$$\left(\int_X(\sigma\bar{\sigma})^n\right)^{n-1}\left(\int_Xc_2Y^{n-1}L^{n-1}\right)^n=\left(\int_XY^nL^n\right)^{n-1}\left(\int_Xc_2(\sigma\bar{\sigma})^{n-1}\right)^n.$$
Note that we have used the fact that $q_X(L)=0$, which
implies that $t_2^{n(n-1)}$ is the highest power of $t_2$
appearing. Next we identify the terms appearing in this equation.

\begin{lem}
\label{first_term}
We have
$$\int_XY^nL^n=n!.$$
\end{lem}

\begin{prf}
Since $L$ is the pullback of a hyperplane in $\P^n$, $L^n$ must be the
pullback of a point, i.e.\ a fibre $F$, which we assume is smooth. The
restriction of $Y$ to $F$ is a theta divisor, and hence
$$\int_XY^nL^n=\int_F\Theta^n=n!$$
since $\Theta$ is a principal polarization of $F$.
\end{prf}

\begin{lem}
\label{second_term}
We have
$$\frac{\left(\int_Xc_2(\sigma\bar{\sigma})^{n-1}\right)^n}{\left(\int_X(\sigma\bar{\sigma})^n\right)^{n-1}}=\frac{24^n(n!)^2}{n^n}\sqrt{\hat{A}}[X]$$
where $\sqrt{\hat{A}}[X]$ is the characteristic number of $X$ coming
from the square root of the $\hat{A}$-polynomial. 
\end{lem}

\begin{rmk}
Note that
\begin{eqnarray*}
\sqrt{\hat{A}} & = & \left(1+\hat{A}_1+\hat{A}_2+\ldots \right)^{1/2} \\
  & = &
  1+\frac{1}{2}\hat{A}_1+\left(\frac{1}{2}\hat{A}_2-\frac{1}{8}\hat{A}_1^2\right)+\ldots \\
  & = & 1+\frac{1}{24}c_2+\frac{1}{5760}(7c_2^2-4c_4)+\ldots . \\
\end{eqnarray*}
In particular $\sqrt{\hat{A}}[X]$ does not mean
$\left(\hat{A}[X]\right)^{1/2}$. 
\end{rmk}

\begin{prf}
The proof of the lemma is based on recognizing that the left hand side
is a {\em Rozansky-Witten invariant} of $X$. Following the notation
of~\cite{hs01}
\begin{eqnarray*}
\int_Xc_2(\sigma\bar{\sigma})^{n-1} & = &
\int_X\frac{1}{16\pi^2n}[\Theta(\Phi)]\sigma^n\bar{\sigma}^{n-1} \\
 & = & \frac{1}{16\pi^2n}c_{\Theta}\int_X(\sigma\bar{\sigma})^n
\end{eqnarray*}
where $\Theta$ denotes the two-vertex trivalent graph and is unrelated
to the theta divisor. Therefore
\begin{eqnarray*}\frac{\left(\int_Xc_2(\sigma\bar{\sigma})^{n-1}\right)^n}{\left(\int_X(\sigma\bar{\sigma})^n\right)^{n-1}}
  & = & \frac{1}{(16\pi^2n)^n}c_{\Theta}^n\int_X(\sigma\bar{\sigma})^n
  \\
  & = & \frac{n!}{2^nn^n}b_{\Theta^n}(X).
\end{eqnarray*}
The main result of Hitchin and the author in~\cite{hs01} is that the
Rozansky-Witten invariant $b_{\Theta^n}(X)$ can be written in terms of
characteristic numbers
$$b_{\Theta^n}(X)=48^nn!\sqrt{\hat{A}}[X]$$
which completes the proof.
\end{prf}

The remaining term $\int_Xc_2Y^{n-1}L^{n-1}$ will be calculated in the
next section.

\section{The second Chern class of $X$}

On $f:X\rightarrow\P^n$ we have the inclusion
$f^*\Omega^1_{\P^n}\rightarrow\Omega^1_X$, which is dual to the
derivative $df:T_X\rightarrow f^*T_{\P^n}$ of $f$. The holomorphic
symplectic form $\sigma$ gives an isomorphism between $\Omega^1_X$ and
$T_X$, so the two maps can be combined into a complex
$$0\rightarrow f^*\Omega^1_{\P^n}\rightarrow\Omega^1_X\cong
T_X\rightarrow f^*T_{\P^n}.$$
For a Lagrangian fibration with good singular fibres, let
$${\mathrm{Sing}}=\cup_{t\in\Delta}\mathrm{Sing}(X_t)$$
be the union of the singular loci of all singular fibres of $X$, and
let $\iota:{\mathrm{Sing}}\hookrightarrow X$ be the inclusion into
$X$. Note that $\mathrm{Sing}$ is a fibration over $\Delta$ whose
generic fibre (over a point of $\Delta_{\mathrm{sm}}$) is an abelian
variety of dimension $n-1$. In particular, $\mathrm{Sing}$ is
codimension two in $X$.

\begin{lem}
\label{Chern}
Let $f:X\rightarrow\P^n$ be a Lagrangian fibration with good singular
fibres. Then
$$0\rightarrow f^*\Omega^1_{\P^n}\rightarrow\Omega^1_X\cong
T_X\rightarrow f^*T_{\P^n}\rightarrow\iota_*{\cal F}\rightarrow 0$$
is exact over $\P^n\backslash\Delta_{\mathrm{sing}}$, where $\cal F$
is a sheaf on $\mathrm{Sing}$ which is generically rank one.
\end{lem}

\begin{prf}
Over smooth fibres and over smooth points of singular fibres our
sequence comes from splicing the two exact sequences
$$\begin{array}{ccccccc}
0 & \rightarrow & f^*\Omega^1_{\P^n} & \rightarrow & \Omega^1_X & 
\rightarrow & \Omega^1_{X/\P^n} \\
 & & & & \downarrow\cong & & \\
0 & \rightarrow & T_{X/\P^n} & \rightarrow & T_X & 
\rightarrow & f^*T_{\P^n}. \\
\end{array}$$
The composition $T_{X/{\P^n}}\rightarrow
T_X\stackrel{\sigma}{\cong}\Omega^1_X\rightarrow\Omega^1_{X/{\P^n}}$
is zero, since $\sigma$ restricted to a (Lagrangian) fibre must
vanish. This proves exactness away from $\mathrm{Sing}$, where all of
the above sheaves are locally free.

In a neighbourhood of $\mathrm{Sing}$ we do a local
computation. Recall that $f$ is given locally by
$$f:(z_1,\ldots,z_n,w_1,\ldots,w_n)\mapsto(z_1w_1,z_2,\ldots,z_n).$$
Therefore $f^*\Omega^1_{\P^n}\rightarrow\Omega^1_X$ is given by
$$\left(\begin{array}{cccccccc}
    w_1 & 0 & \ldots & 0 \\
    0   & 1 & \ldots & 0 \\
    \vdots & \vdots & \ddots & \vdots \\
    0   & 0 & \ldots & 1 \\
    z_1 & 0 & \ldots & 0 \\
    0  & 0 & \ldots & 0 \\
    \vdots & \vdots & \ddots & \vdots \\
    0  & 0 & \ldots & 0 \\
\end{array}\right),$$
the isomorphism $\Omega^1_X\stackrel{\sigma}{\cong}T_X$ is given by
$$\left(\begin{array}{cc}
    0 & \mathrm{Id}_{n\times n} \\
    -\mathrm{Id}_{n\times n} & 0 \\
\end{array}\right),$$
and $df:T_X\rightarrow f^*T_{\P^n}$ is given by
$$df=\left(\begin{array}{cccccccc}
    w_1 & 0 & \ldots & 0 & z_1 & 0 & \ldots & 0 \\
    0   & 1 & \ldots & 0 & 0   & 0 & \ldots & 0 \\
    \vdots & \vdots & \ddots & \vdots & \vdots & \vdots & \ddots & \vdots \\
    0   & 0 & \ldots & 1 & 0   & 0 & \ldots & 0 \\
\end{array}\right).$$
It is now a simple matter to check that
$$0\rightarrow f^*\Omega^1_{\P^n}\rightarrow\Omega^1_X\cong
T_X\rightarrow f^*T_{\P^n}$$
is exact, and that $df$ drops rank by one when $w_1=z_1=0$, which are
precisely the local equations for $\mathrm{Sing}$. Thus the cokernel
of $df$ looks like $\iota_*{\cal F}$ where $\cal F$ is a generically
rank one sheaf on $\mathrm{Sing}$.
\end{prf}

It follows immediately from the lemma that
$$c_1(T_X)=f^*c_1(\Omega^1_{\P^n})+f^*c_1(T_{\P^n})=0$$
and
\begin{eqnarray*}
c_2(T_X) & = & [{\mathrm{Sing}}]+f^*c_2(\Omega^1_{\P^n})+f^*c_2(T_{\P^n})+\mbox{const.}[f^{-1}(\Delta_{\mathrm{sing}})] \\
   & = & [{\mathrm{Sing}}]+n(n+1)L^2+\mbox{const.}L^2 
\end{eqnarray*}
for some constant.

\begin{rmk}
This formula for the second Chern class is the holomorphic analogue of
a well-known formula relating the first Chern class and singular locus
of a real Lagrangian fibration on a (real) symplectic manifold. It is
really the key to Theorem~\ref{principal} below, as $[\mathrm{Sing}]$
will lead directly to $\mathrm{deg}\Delta$, while we already saw that
$c_2(T_X)$ leads to $\sqrt{\hat{A}}[X]$. 
\end{rmk}

\begin{lem}
\label{third_term}
We have
$$\int_Xc_2Y^{n-1}L^{n-1}=\int_X[{\mathrm{Sing}}]Y^{n-1}L^{n-1}=(n-1)!\mathrm{deg}\Delta.$$
\end{lem}

\begin{prf}
Firstly
\begin{eqnarray*}
\int_Xc_2Y^{n-1}L^{n-1} & = & \int_X[{\mathrm{Sing}}]Y^{n-1}L^{n-1}+\mbox{const.}Y^{n-1}L^{n+1} \\
 & = & \int_X[{\mathrm{Sing}}]Y^{n-1}L^{n-1}
\end{eqnarray*}
since $L^{n+1}=0$ ($L$ is the pull-back of a divisor from the
$n$-dimensional base).

The locus $\mathrm{Sing}$ is supported over the discriminant locus
$\Delta$, while $L^{n-1}$ is the pull-back of a line $\ell$ in
$\P^n$. Since we can assume $\ell$ is generic, it will intersect
$\Delta$ in precisely $\mathrm{deg}\Delta$ points, with each point in
$\Delta_{\mathrm{sm}}$. In this way we reduce the lemma to computing
an intersection number in a good singular fibre. This computation
will be invariant under deformation, so we can assume that the good
singular fibre is the compactified Jacobian $\bar{J}_0C$ of a curve
$C$ with a single node. 

The restriction of $\mathrm{Sing}$ to $\bar{J}_0C$ is of course the
singular locus $s$ which comes from identifying $s_0$ and
$s_{\infty}$. The restriction of $Y$ to $\bar{J}_0C$ is the
generalized theta divisor $\Theta$. In the Jacobian $J_0C$ of a smooth
curve $C$, $\Theta^{n-1}$ is cohomologous to $(n-1)!C$, with $C$
embedded in $J_0C$ by the Abel-Jacobi map. In fact this relation
remains true for a curve with a single node, which can also be
embedded in its compactified Jacobian by a generalization of the
Abel-Jacobi map. Then $C$ intersects the singular locus $s$ at
precisely one point, the node of $C$.

Combining the above observations we find
$$\int_X[{\mathrm{Sing}}]Y^{n-1}L^{n-1}=\mathrm{deg}\Delta\int_{\bar{J}_0C}[s]\Theta^{n-1}=(n-1)!\mathrm{deg}\Delta.$$
\end{prf}

\begin{rmk}
One could also observe that the restriction of the generalized theta
divisor $\Theta$ to the singular locus $s$ induces a principal
polarization on $s$, and thus 
$$\int_{X_t}[s]\Theta^{n-1}=\int_s(\Theta|_s)^{n-1}=(n-1)!.$$
There is then no need to mention compactified Jacobians.
\end{rmk}

These calculations now yield a formula for the degree of $\Delta$.
\begin{thm}
\label{principal}
Let $X\rightarrow\P^n$ be a Lagrangian fibration by principally
polarized abelian varieties, by which we mean that there is a divisor
$Y$ on $X$ which restricts to (a multiple of) a principal polarization
on the generic fibre. If $X$ has good singular fibres then
\begin{eqnarray*}
\mathrm{deg}\Delta & = & \frac{1}{2}b_{\Theta^n}(X)^{\frac{1}{n}} \\
 & = & 24\left(n!\sqrt{\hat{A}}[X]\right)^{\frac{1}{n}}.
\end{eqnarray*}
\end{thm}

\begin{prf}
We simply substitute the results of Lemmas~\ref{first_term},
\ref{second_term}, and \ref{third_term} into the equation preceding
Lemma~\ref{first_term}. Note that even if $Y$ restricts to a
non-trivial multiple $m\Theta$ of a theta divisor on each fibre, the
factor $m$ will ultimately cancel out.
\end{prf}

\begin{rmk}
The hypotheses imply that $X$ is projective, as $Y+kL$ will be ample
for sufficiently large $k$. However, we expect that the formula will
hold more generally, when the generic fibre is only abstractly a
principally polarized abelian variety, without any reference to a
global divisor on $X$. The reason is that there are ways to deform a
Lagrangian fibration until it admits a section (see~\cite{sawon04}
and~\cite{sawon05}) without changing the local structure of the
fibration, and in particular, without changing the discriminant locus
$\Delta$. Now a Lagrangian fibration is projective if and only if it
admits a rational section or multi-section (Proposition 5.1 of
Oguiso~\cite{oguiso06}). In particular, our Lagrangian fibration with
a section will contain an ample divisor $Y$, which should then induce
the principal polarization of the generic fibre.
\end{rmk}

\section{The Beauville-Mukai system}

In this section we verify our formula for the Beauville-Mukai
integrable system~\cite{beauville99} described in Section 2, whose
total space is a deformation of the Hilbert scheme $S^{[n]}$ of $n$
points on a K3 surface $S$. In~\cite{sawon00} the author calculated
various Rozansky-Witten invariants; in particular
$$b_{\Theta^n}(S^{[n]})=12^n(n+3)^n.$$
Applying Theorem~\ref{principal}, the discriminant locus of a
fibration on $S^{[n]}$ (or on any deformation of $S^{[n]}$) should
therefore have degree
$$\mathrm{deg}\Delta=6(n+3).$$

For $n=1$ it is well-known that a generic elliptic K3 surface has
exactly $24$ singular fibres. For $n\geq 2$ we have the
Beauville-Mukai system coming from a genus $n$ curve $C$ contained in
$S$, which is a fibration over $|C|\cong\P^n$. There is a map
$S\rightarrow(\P^n)^{\vee}$ which for generic $S$ is an embedding (or
branched double cover when $n=2$). The discriminant locus
$\Delta\subset |C|$ parametrizes singular curves in the linear system,
i.e.\ it parametrizes hyperplanes in $(\P^n)^{\vee}$ whose
intersection with $S$ is singular. In other words, $\Delta\subset\P^n$
is the variety dual to $S\subset(\P^n)^{\vee}$ (or dual to the branch
curve of $S\rightarrow(\P^2)^{\vee}$ when $n=2$).  

Consider a pencil of
hyperplanes $H_t\subset(\P^n)^{\vee}$, with $t\in\P^1$. Generically
there will be $\mathrm{deg}\Delta$ singular hyperplane sections of $S$
in this pencil, and each one will have a single node. The union
$\cup_{t\in\P^1}H_t\cap S$ of these hyperplane sections
gives a divisor in $S\times\P^1$ whose corresponding line bundle is
${\cal O}(C,1)$. If this divisor is given locally by $f=0$, then the
singularities of $H_t\cap S$ are given by $f=0$ and $df=0$, where the
derivative is taken only in the direction of $S$. Globally, we have a
section of the rank three vector bundle
$${\cal O}(C,1)\oplus T^*S(C,1)$$
which vanishes precisely at the singular points. Therefore
\begin{eqnarray*}
\mathrm{deg}\Delta & = & c_3({\cal O}(C,1)\oplus
T^*S(C,1))[S\times\P^1] \\
 & = & 6(n+3)
\end{eqnarray*}
where we have used the fact that $C^2=2n-2$. Thus we have a
verification of Theorem~\ref{principal} in this case.

\section{Non-principal polarizations}

Let us illustrate how to modify our theorem for non-principal
polarizations. Let $X$ be an irreducible holomorphic symplectic
manifold fibred over $\P^n$, and let $Y$ be a divisor on $X$ which
on a generic fibre restricts to a polarization of type
$(d_1,\ldots,d_n)$ with $d_1|d_2|\cdots|d_n$. We first generalize the
notion of a good singular fibre to this case: in fact there is more
than one model.


A singular fibre $X_t$, with $t\in\Delta_{\mathrm{sm}}$, should look
like a generic semi-stable degeneration of an abelian variety with
polarization of type $(d_1,\ldots,d_n)$. In other words, $X_t$ should
be a semi-stable degeneration that occurs over a generic point of the
boundary of an Igusa/Mumford~\cite{igusa67,mumford72,mumford75}
compactification of the moduli space of abelian varieties. For
non-principal polarizations, the boundary consists of several
irreducible (codimension one) components, thus we expect to find
several different models which we now describe explicitly.

The normalization $\tilde{X_t}$ of $X_t$ will look like a collection
of $k$ $\P^1$-bundles over an abelian variety of dimension $n-1$. The
singular fibre itself is obtained by gluing the zero and infinity
sections in a chain, as shown (with $s^1_0$ also glued to
$s^k_{\infty}$, with a translation). 
$$\begin{array}{cccl}
 & & & s^k_{\infty} \\
\mathbb{P}^1 & \hookrightarrow & \singularI & \\
 & & & s^{k-1}_{\infty}\cong s^k_0 \\
 & & \vdots & \\
 & & \vdots\hspace*{20mm}\vdots & \\
 & & \vdots & \\
 & & & s^1_{\infty}\cong s^2_0 \\
\mathbb{P}^1 & \hookrightarrow & \singularII \\
 & & & s^1_0
\end{array}$$
Note that the singular locus $\mathrm{Sing}(X_t)$ consists of $k$
irreducible components, each isomorphic to the abelian variety of
dimension $n-1$. Moreover the polarization of a nearby smooth fibre,
which is of type $(d_1,\ldots,d_n)$, will degenerates to a divisor
$Y_t$ in $X_t$. Suppose that $Y_t$ induces a polarization of type
$(d^{\prime}_1,\ldots,d^{\prime}_{n-1})$ on each irreducible component
of $\mathrm{Sing}(X_t)$. Compatibility requires that
$d_i|d^{\prime}_i$ for $i=1,\ldots,n-1$, and
$$d_1d_2\cdots d_{n-1}d_n=d^{\prime}_1d^{\prime}_2\cdots
d^{\prime}_{n-1}k.$$
In particular, this implies that $k$ must divide $d_n$. For example,
in the case of abelian surfaces with polarization of type $(1,p)$,
with $p$ prime, there are two possible degenerations: one is
irreducible whereas the other consists of $p$ irreducible components
(see Propositions 4.5 and 4.7 in Hulek, Kahn, and
Weintraub~\cite{hkw93}).

\begin{dfn}
We say a Lagrangian fibration $X\rightarrow\P^n$ by abelian varieties
with polarization of type $(d_1,\ldots,d_n)$ has {\em good} singular
fibres if the generic singular fibre $X_t$ for
$t\in\Delta_{\mathrm{sm}}$ looks like the picture described
above. Note that $\Delta$ may consist of several irreducible
components and the model for the generic singular fibre $X_t$ may
differ over each component (e.g.\ $k$ and
$(d^{\prime}_1,\ldots,d^{\prime}_{n-1})$ need not be the same over
every component).
\end{dfn}

Let $L$ be the pullback of a hyperplane in $\P^n$, and $Y$ the
relative theta divisor. As before, we have
$$\left(\int_X(\sigma\bar{\sigma})^n\right)^{n-1}\left(\int_Xc_2Y^{n-1}L^{n-1}\right)^n=\left(\int_XY^nL^n\right)^{n-1}\left(\int_Xc_2(\sigma\bar{\sigma})^{n-1}\right)^n.$$
Lemma~\ref{first_term} becomes
$$\int_XY^nL^n=n!d_1d_2\cdots d_{n-1}d_n.$$
Lemma~\ref{second_term} remains unchanged
$$\frac{\left(\int_Xc_2(\sigma\bar{\sigma})^{n-1}\right)^n}{\left(\int_X(\sigma\bar{\sigma})^n\right)^{n-1}}=\frac{24^n(n!)^2}{n^n}\sqrt{\hat{A}}[X].$$
The exact sequence of Lemma~\ref{Chern} also remains unchanged,
because although the singular locus $\mathrm{Sing}(X_t)$ of each
generic singular fibre now consists of $k$ irreducible components, the
local description of $f:X\rightarrow\P^n$ near these singularities
does not change. Therefore our expression for the second Chern class
of $X$ is still valid, and Lemma~\ref{third_term} becomes
\begin{eqnarray*}
\int_Xc_2Y^{n-1}L^{n-1} & = & \int_X[\mathrm{Sing}]Y^{n-1}L^{n-1} \\
 & = & \mathrm{deg}\Delta\int_{X_t}[\mathrm{Sing}(X_t)]Y^{n-1}_t \\
 & = & k(n-1)!d^{\prime}_1d^{\prime}_2\cdots
 d^{\prime}_{n-1}\mathrm{deg}\Delta \\
 & = & (n-1)!d_1d_2\cdots d_n\mathrm{deg}\Delta \\
\end{eqnarray*}
because $\mathrm{Sing}(X_t)$ consists of $k$ irreducible components,
each isomorphic to an abelian variety of dimension $n-1$, and $Y_t$
intersects each component in a polarization of type
$(d^{\prime}_1,\ldots,d^{\prime}_{n-1})$.

Combining these formulae we obtain the following result.
\begin{thm}
\label{nonprincipal}
Let $X\rightarrow\P^n$ be a Lagrangian fibration by abelian varieties
with polarization of type $(d_1,\ldots,d_n)$, by which we mean that
there is a divisor $Y$ on $X$ which restricts to a polarization of
this type on the generic fibre. If $X$ has good singular fibres then 
\begin{eqnarray*}
\mathrm{deg}\Delta & = &
\frac{1}{2}\left(\frac{b_{\Theta^n}(X)}{d_1\cdots d_n}\right)^{\frac{1}{n}} \\
 & = & 24\left(\frac{n!\sqrt{\hat{A}}[X]}{d_1\cdots
     d_n}\right)^{\frac{1}{n}}.
\end{eqnarray*}
\end{thm}

\begin{rmk}
One could always change the polarization of $X$ to $mY$ with $m\geq
2$, and this would multiply all the $d_i$ by the factor $m$. Our
formula then appears to be inconsistent; however, our models for
singular fibres
implicitly assume that $Y$ is a primitive divisor. This suggests that
we should assume $d_1=1$. Indeed if $d_1>1$ then $Y$ is not primitive
when restricted to a fibre, and in some circumstances one can use the
methods described in~\cite{sawon04} and~\cite{sawon05} to deform
$X\rightarrow\P^n$ so that $Y=d_1Y^{\prime}$ globally, without
changing the fibration locally. Changing to the new polarization
$Y^{\prime}$, we could then assume that $d_1=1$.
\end{rmk}

\section{Generalized Kummer varieties}

The generalized Kummer varieties $K_n$ were introduced by
Beauville~\cite{beauville83}. Debarre~\cite{debarre99} exhibited
a fibration on $K_n$; see also Example 3.8 in~\cite{sawon03}. The
fibres have polarization of type $(1,\ldots,1,n+1)$ and this fibration
has good singular fibres. In~\cite{sawon00} the author calculated
$$b_{\Theta^n}(K_n)=12^n(n+1)^{n+1}.$$
Theorem~\ref{nonprincipal} therefore gives
$$\mathrm{deg}\Delta=6(n+1).$$

For $n=1$ this gives twelve. This is correct because the Kummer K3
surface $K_1$ will be an elliptic fibration whose singular fibres each
consist of two irreducible components; more precisely, they are of
Kodaira type $I_2$ and so there will indeed be twelve of them.

For $n\geq 2$ one begins with an abelian surface $A$ with
polarization of type $(1,n+1)$. Thus $A$ is polarized by a genus
$n+2$ curve $C$ with $C^2=2(n+1)$. The relative Jacobian of the family
of curves linear equivalent to $C$ is a fibration over $|C|\cong\P^n$
whose generic fibre is an abelian variety of dimension $n+2$. There is
a map from the total space of this fibration to $A$ (the Albanese
map), and the kernel of this map gives a fibration on $K_n$. More
precisely, the kernel is isomorphic to the generalized Kummer variety
$K_n(\hat{A})$ constructed from the dual abelian surface $\hat{A}$,
and it inherits the map to $\P^n$ which makes it a Lagrangian
fibration.

As with the Beauville-Mukai system, $\Delta\subset\P^n$ parametrizes
hyperplanes in $(\P^n)^{\vee}$ whose intersection with $A\subset
|C|^{\vee}\cong(\P^n)^{\vee}$ is singular (with the obvious
modifications for small $n$, when $A$ is not necessarily embedded). We
can therefore use the same method to calculate the degree of $\Delta$,
and we obtain
\begin{eqnarray*}
\mathrm{deg}\Delta & = & c_3({\cal O}(C,1)\oplus
T^*A(C,1))[A\times\P^1] \\
 & = & 6(n+1)
\end{eqnarray*}
which agrees with the value obtained from Theorem~\ref{nonprincipal}.


\section{Fibrations on four-folds}

Suppose $X\rightarrow\P^2$ is an irreducible holomorphic symplectic
four-fold which admits a Lagrangian fibration by abelian surfaces with
polarization of type $(d_1,d_2)$, and write $d_2=d_1d$. Moreover,
let's follow the remark after Theorem~\ref{nonprincipal} and assume
$d_1=1$. In this dimension, if the base is smooth then Matsushita's
results imply it must be $\P^2$. If the fibration has good singular
fibres then Theorem~\ref{nonprincipal} yields
$$\mathrm{deg}\Delta=\frac{1}{2}\left(\frac{b_{\Theta^n}(X)}{d}\right)^{\frac{1}{2}}=\left(\frac{1152\sqrt{\hat{A}}[X]}{d}\right)^{\frac{1}{2}}.$$
We will use Guan's bounds~\cite{guan01} on the Betti numbers of
$X$ to restrict the possible values of $d$ and $\mathrm{deg}\Delta$.

\begin{thm}[Guan~\cite{guan01}]
Let $X$ be an irreducible holomorphic symplectic four-fold. The Betti
numbers of $X$ are bounded and can only take the following values:
\begin{itemize}
\item $b_2=23$ and $b_3=0$,
\item $b_2=8$ and $b_3=0$,
\item $b_2=7$ and $b_3=0$ or $8$,
\item $b_2=6$ and $b_3=0$, $4$, $8$, $12$, or $16$,
\item $b_2=5$ and $b_3=0$, $4$, $8$, $\ldots$ or $36$,
\item $b_2=4$ and $b_3=0$, $4$, $8$, $\ldots$ or $60$,
\item $b_2=3$ and $b_3=0$, $4$, $8$, $\ldots$ or $68$.
\end{itemize}
\end{thm}
The fourth Betti number is determined by Salamon's relation
$$b_4=46+10b_2-b_3$$
and therefore
$$c_4[X]=\chi(X)=48+12b_2-3b_3.$$
The relation
$$\hat{A}[X]=\frac{1}{720}(3c_2^2[X]-c_4[X])=\chi(\O_X)=3$$
between the Chern numbers allows us to write $\sqrt{\hat{A}}[X]$
solely in terms of $c_4[X]$, giving
$$1152\sqrt{\hat{A}}[X]=1008-\frac{1}{3}c_4[X]=992-4b_2+b_3.$$
We can now state our final result.
\begin{thm}
\label{values}
Let $X$ be an irreducible holomorphic symplectic four-fold which
admits a Lagrangian fibration $X\rightarrow\P^2$ by abelian surfaces
with polarization of type $(1,d)$, by which we mean that there is a
divisor $Y$ on $X$ which restricts to a polarization of type $(1,d)$
on the generic fibre. If $X$ has good singular fibres then
$\mathrm{deg}\Delta$ is at most $32$ and $d$ is at most $1036$.
\end{thm}

\begin{prf}
Firstly, $b_2$ must be at least four since $L$ corresponds to an
isotropic element of $\H^2(X,\mathbb{Z})$ with respect to the
Beauville-Bogomolov form, and this is a lattice of signature
$(3,b_2-3)$. We substitute the possible values of $b_2$ and $b_3$ (as
allowed by Guan's Theorem) into
$$1152\sqrt{\hat{A}}[X]=992-4b_2+b_3.$$
The largest value is $1036$ when $b_2=4$ and $b_3=60$. Moreover,
$1152\sqrt{\hat{A}}[X]$ is always an integer and the formula for
$\mathrm{deg}\Delta$ shows that it must be divisible by $d$. Thus $d$
is at most $1036$. Moreover
$$\mathrm{deg}\Delta=\left(\frac{1152\sqrt{\hat{A}}[X]}{d}\right)^{\frac{1}{2}}\leq
\left(\frac{1036}{d}\right)^{\frac{1}{2}}\leq\sqrt{1036}<33.$$
\end{prf}

\begin{rmk}
In our two examples we have $d=1$ and $\mathrm{deg}\Delta=30$ for the
Beauville-Mukai system on $S^{[2]}$, and $d=3$ and
$\mathrm{deg}\Delta=18$ for the generalized Kummer four-fold $K_2$. We
suspect that further work will eliminate many (perhaps all) of the
other possible values for $d$ and $\mathrm{deg}\Delta$. In particular,
it is hard to imagine there could be any examples with $d$ large.
\end{rmk}

\begin{flushleft}
Department of Mathematics\hfill sawon@math.sunysb.edu\\
SUNY at Stony Brook\hfill www.math.sunysb.edu/$\sim$sawon\\
Stony Brook NY 11794-3651\\
USA\\
\end{flushleft}

\end{document}